\documentclass[12pt]{tran-l}
\usepackage{diagrams}

\diagramstyle[height=2em,width=2em]
\newarrow{Epi}----{>>}
\newarrow{Mono}>--->
\newarrow{Iso}>---{>>}
\newarrow{Mapsto}|--->
\newarrow{Igual}=====
\newarrow{Dashto}{}{dash}{}{dash}>

\newcommand{\Z}{{\mathbb Z}}

\renewcommand{\H}{{\mathcal H}}
\newcommand{\F}{{\mathbb F}}

\newcommand{\pcom}{_{p}^{\wedge}}
\newcommand{\doscom}{_{2}^{\wedge}}

 \renewcommand{\P}[1]{\mathcal{P}^{#1}}

\newcommand{\hocolim}{\operatornamewithlimits{hocolim}}

\newcommand{\map}{\operatorname{Map}\nolimits}

\newcommand{\A}{\ifmmode{\mathcal{A}}\else${\mathcal{A}}$\fi}
\newcommand{\K}{\ifmmode{\mathcal{K}}\else${\mathcal{K}}$\fi}
\newcommand{\U}{\ifmmode{\mathcal{U}}\else${\mathcal{U}}$\fi}
\newcommand{\Qp}{\ifmmode{\mathcal{Q}}\else${\mathcal{Q}}$\fi}
\newtheorem{Thm}{Theorem}[section]
\newtheorem{Prop}[Thm]{Proposition}

\newtheorem{Lem}[Thm]{Lemma}

\newtheorem{Con}[Thm]{Conjecture}

\theoremstyle{definition}

\newtheorem{Rmk}[Thm]{Remark}

\newtheorem*{teoA}{Theorem A}
\newtheorem*{teoB}{Theorem B}
\newtheorem*{teoC}{Theorem C}

\theoremstyle{remark}

\newcommand{\diag}{\operatorname{Diag}\nolimits}

\newcommand{\Mor}{\operatorname{Mor}\nolimits}

\title[On the mod~$p$ cohomology of $BPU(p)$]
{On the mod~$p$ cohomology of $BPU(p)$}

\thanks{The first
author is partially supported by the Ministry for Education, Science and Sport of
Republic of Slovenia Research Program no.~101-509. The second author is partially
supported by the DGES-FEDER grant BFM2001-1825, and Junta de Andaluc{\'\i}a Grant
FQM-0213.}

\author{Ale\v s Vavpeti\v c} \address{\hfill\break Fakulteta za matematiko in fiziko\\
\hfill\break Univerza v Ljubljani\\ \hfill\break Jadranska 19\\ \hfill\break SI-1111
Ljubljana\\ \hfill\break Slovenia} \email{ales.vavpetic@FMF.Uni-Lj.Si}
\author{Antonio Viruel}
\address{\hfill\break Dpto de {\'A}lgebra, Geometr{\'\i}a y Topolog{\'\i}a\\
\hfill\break Universidad de M{\'a}laga\\ \hfill\break Apdo correos 59 \\ \hfill\break
E29080 M{\'a}laga
\\ \hfill\break Spain} \email{viruel@agt.cie.uma.es}

\date{\today}

\subjclass[2000]{55R35, 55R15}

\begin{document}

\begin{abstract} We study the mod~$p$ cohomology of the classifying
space of the projective unitary group $PU(p)$. We first proof that old conjectures due
to J.F.~Adams, and Kono and Yagita \cite{KY} about the structure of the mod~$p$
cohomology of classifying space of connected compact Lie groups held in the case of
$PU(p)$. Finally, we proof that the classifying space of the projective unitary group
$PU(p)$ is determined by its mod~$p$ cohomology as an unstable algebra over the
Steenrod algebra for $p>3$, completing previous works \cite{DMW1} and \cite{BV} for
the cases $p=2,3$.
\end{abstract}

\maketitle

\section{Introduction}

Compact Lie groups provide an example of one the classical mathematical maxims: ``the
richer is the mathematical structure of an object, the more rigid the object is". So
for example all the rich mathematical structure associated to a connected compact Lie
group is so intimately linked that it is completely recover (perhaps, up to local
isomorphism) from some small data like the Dynkin diagram or the maximal torus
normalizer \cite{CWW}.

In homotopy theory, this rigidity in the structure of a compact Lie group $G$ is
expected to be inherited by the classifying space $BG$ and related structures. So for
example, in the appropriate homotopical setting of $p$-compact groups \cite{DW},
maximal torus normalizers do characterize the isomorphic type of $BG$, at least at odd
primes \cite{AGMV}.

The aim of this work is to study the rigidity of the mod~$p$
cohomology of $BG$, namely $H^*(BG;\F_p)$, proving several
conjectures in the particular case of $G$ being the projective
unitary group $PU(p)$, which is obtained as the quotient of the
unitary group of rank $p$, $U(p)$, by the subgroup $\{
\diag(\alpha,\ldots,\alpha)\mid \alpha \in S\sp 1\}$ of diagonal
matrices.

In \cite[Theorem 1.1]{DW}, it is shown that $H^*(BG;\F_p)$ is a Noetherian algebra for
any compact connected Lie group $G$, so by \cite[Theorem 1.4]{R} (or directly
\cite[Theorem 6.2]{Q1}) we know that the kernel of the natural map
\begin{equation}\label{adams}H^*(BG;\F_p)\rTo \lim_{\overleftarrow{\A_p(G)}}
H^*(BE;\F_p),\end{equation} where $\A_p(G)$ stands for the Quillen category of
elementary abelian $p$-subgroups of $G$ \cite{Q1,R,JM,DW-decomposition}, contains only
nilpotent elements. For $p>2$, a more stronger conjecture shows up

\begin{Con}[J.F.\ Adams]\label{adams-conj}
Let $G$ be a compact connected Lie group, and $p$ be an odd prime.
Then the mod~$p$ cohomology of $BG$ is detected by elementary
abelian $p$-subgroups \cite[Definition 4.2]{AM}, i.e.\ the natural
map \eqref{adams} is a monomorphism.
\end{Con}

Conjecture \ref{adams-conj} trivially holds in the torsion free cases (see
\cite[Theorem 12.1]{AGMV}). In the case of torsion, only a few examples have been
worked out, all of them at $p=3$: $F_4$ \cite[Teorema 5]{Broto}, $E_6$ \cite{mstt} and
$PU(3)$ \cite[Theorem 3.3]{KY}. Our first result generalizes the last reference
proving,

\begin{teoA}
The group $PU(p)$ verifies Conjecture \ref{adams-conj} at the odd prime $p$, i.e.
$H^*(BPU(p);\F_p)$ is detected by elementary abelian $p$-subgroups, i.e.\ the natural
map \eqref{adams} is a monomorphism for the case $G=PU(p)$.
\end{teoA}
\begin{proof}
See Theorem \ref{demo-adams}
\end{proof}

The knowledge of the structure of $H^*(BG;\F_p)$ plays an important role when trying
to understand other generalized cohomologies of $BG$ as it is shown in \cite{KY}. So
for example, understanding Milnor primitive operations (see Section \ref{KY-sec}) is a
crucial step in the use of the Atiyah-Hirzebruch spectral sequence
\cite[pag.~496]{McCleary}. A new conjecture arises \cite[Conjecture 5]{KY}

\begin{Con}[Kono-Yagita]\label{KY-conj}
If $G$ is a connected compact Lie group, then for each odd dimensional element $x\in
H^*(BG;\F_p)$, there is $i$ such that such that $Q_mx\ne 0$ for all $m\ge i$, where
$Q_m$ are the Milnor primitive operators.
\end{Con}

Then our result generalizes the study of $PU(3)$ carried out in \cite{KY} proving,

\begin{teoB}
The group $PU(p)$ verifies Conjecture \ref{KY-conj} for every odd prime $p$, i.e. for
each odd dimensional element $x\in H^*(BPU(p);\F_p)$, there is $i$ such that 
$Q_mx\ne 0$ for all $m\ge i$, where $Q_m$ are the Milnor primitive operators.
\end{teoB}
\begin{proof}
See Theorem \ref{demo-KY}.
\end{proof}

\begin{Rmk}
It is worth to remark that while the proof of Conjectures \ref{adams-conj} and
\ref{KY-conj} in previous known cases is heavily based in a precise understanding of
the cohomology rings involved, i.e.\ generators and relations, the proofs of Theorems
A and B is done by geometrical methods and without using any information about the
algebra structure of $H^*BPU(p)$.
\end{Rmk}

So many restrictions on $H^*(BG;\F_p)$ suggest that these algebras
do not show up in nature very frequently. In other words, any
space $X$ whose mod~$p$ cohomology is isomorphic to that of $BG$,
for a connected compact Lie group $G$, should be topologically
related with $BG$ in some way. This idea is captured in the next
conjecture \cite[Conjecture 4.4]{N-survey}

\begin{Con}\label{notbohm}
Let $G$ be a compact connected Lie group, and let $X$ be a $p$-complete space such
that $H^*(X;\F_p)\cong H^*(BG;\F_p)$ as algebras over the mod~$p$ Steenrod algebra
$\A_p$. Then $X\simeq BG\pcom$.
\end{Con}

The first result of this kind appeared in \cite{DMW}, where Dwyer, Miller and
Wilkerson proved Conjecture \ref{notbohm} for $G=SU(2)=S^3$ at $p=2$. In \cite{DMW1},
the same authors considered the case when $p$ does not divide the order of the Weyl
group of $G$. Notbohm in \cite{N} considered the case when $p$ divides the order of
the Weyl group of $G$, but $BG$ has no torsion. For the case when torsion exists,
there are only few known results \cite{BV,V,V2,V3}. We prove,

\begin{teoC}
Let $X$ be a $p$-complete space such that $H^*(X;\F_p)\cong H^*(BPU(p);\F_p)$ as an
unstable algebra over the Steenrod algebra $\mathcal{A}_p$. Then $X$ is homotopy
equivalent to $BPU(p)\pcom$.
\end{teoC}
\begin{proof}
If $p=2$, then $PU(2)=SO(3)$, and the theorem is known \cite{DMW}. If $p=3$ the
theorem is proved in \cite{BV}. The case $p\ge 5$ is consider in Section \ref{demo-N}.
\end{proof}

\noindent{\bf Notation:} Here all spaces are assumed to have the homotopy type of
CW-complexes. Completion means Bousfield-Kan completion \cite{BK}. For a given space
$X$, we write $H^*X$ for the mod~$p$ cohomology $H^*(X;\F_p)$ and $X\pcom$ for
Bousfield-Kan $(\Z_p)_\infty$-completion or $p$-completion of the space $X$. Given a
group $G$ and a $\Z G$-module $M$, we write $\H^*(G;M)$ for the cohomology of $G$ with
(twisted) coefficients in $M$. We assume that the reader is familiar with the Lannes'
theory \cite{L}.


\section{The Adams' conjecture}

The aim of this section is to prove Adams' conjecture (Conjecture \ref{adams-conj})
for the group $PU(p)$ at the prime $p>2$. We start identifying some important
subgroups of a compact connected Lie group $G$. Let $T(G)\subset G$ be a
maximal torus and $N(G)\subset G$ its normalizer. Define $N_p(G)\subset N(G)$, the
$p$-normalizer of the maximal torus $T(G)$, as the preimage of a $p$-Sylow subgroup in
the Weyl group of $G$, $W_G=N(G)/T(G)$.

\begin{Lem}\label{l1}
The groups $N_p(PU(p))$ and $N_p(SU(p))$ are isomorphic.
\end{Lem}
\begin{proof}
Notice first that $N_p(PU(p))=N_p(SU(p))/\{ \diag(\alpha,\ldots,\alpha)\mid \alpha \in
S\sp 1\}$. Now, every element in $N_p(SU(p))$ can be written in a unique way as
$\diag(z_1,\ldots,z_p)P^i$, where $P$ is the permutation matrix corresponding to the
cyclic permutation $(1,2,\ldots,p)$. Then $\varphi\colon N_p(PU(p))\rTo N_p(SU(p))$
given by $$ \varphi([\diag(z_1,\ldots,z_p)P^i])=\diag({z_1\over
z_2},\ldots,{z_{p-1}\over z_p},{z_p\over z_1})P^i $$ provides the desired isomorphism.
\end{proof}

The isomorphism constructed above is very convenient as we can prove Conjecture
\ref{adams-conj} for $N_p(SU(p))$,

\begin{Lem}\label{l2}
The cohomology $H^*BN_p(SU(p))$ is detected by elementary abelian subgroups.
\end{Lem}
\begin{proof}
Notice tha $N_p(U(p))\cong S^1\wr \Z/p$, hence by \cite[Theorem 4.3]{AM} we know that
$H^*BN_p(U(p))$ is detected by elementary abelian subgroups. Moreover $H^*BN_p(U(p))$
is detected by just two subgroups, $V_t=\Z/p^{\oplus p}\subset T(U(p))$ the maximal
elementary abelian toral subgroup and
$V_n=\Z/p\oplus\Z/p=\big(Z(U(p))\times\Z/p\big)\cap SU(p)$ \cite[Lemma 4.4]{AM},
where $Z(U(p))$ is the center of the group $U(p)$

Now the fibration $$S^1\rTo BSU(p)\rTo BU(p)$$ gives rise to a fibration $$ S^1\rTo
BN_p(SU(p))\rTo^{Bj} BN_p(U(p)),$$ whose the Gysin sequence is $$ \cdots\rTo
H^*BN_p(U(p))\rTo^{Bj^*} H^*BN_p(SU(p))\rTo^{d} H^{*-1}BN_p(U(p))\rTo\cdots $$ Let
$x\in H^*BN_p(SU(p))$ and suppose $d(x)\ne 0$. Let $V\rTo BN_p(U(p))$ be an elementary
abelian group detecting $d(x)$. Then $V'=\langle V, Z(U(p))\rangle \cap BN_p(SU(p))$
is an elementary abelian group, which appears in the fibration 
$S^1\rTo BV'\rTo B\langle V, Z(U(p))\rangle$, and detects the element $x$.

If $d(x)=0$, then $x=Bj^*(y)$ for some $y\in H^*BN_p(U(p))$ and $y$ is detected by
$V_t$ or $V_n$ defined above, so the element $x$ is detected by $V_n$ or $V_t\cap
N_p(SU(p))\cong (\Z/p)^{p-1}$.
\end{proof}

An easy consequence of the previous lemmas is

\begin{Lem}\label{DetElAb}
The mod~$p$ cohomology of $BN(PU(p))$ is detected by elementary abelian $p$-subgroups.
\end{Lem}
\begin{proof}
Combining Lemmas \ref{l1} and \ref{l2} we obtain that $H^*BN_p(PU(p))$ is detected by
elementary abelian $p$-subgroups. Then, because the index $[N(PU(p))\colon
N_p(PU(p))]=(p-1)!$ is nonzero in $\F_p$, the transfer argument \cite[Lemma
6.7.17]{WE} shows that $H^*BN(PU(p))\rTo H^*BN_p(PU(p))$ is a monomorphism. Therefore
$H^*BN(PU(p))$ is also detected by elementary abelian $p$-subgroups.
\end{proof}

Finally,

\begin{Thm}\label{demo-adams}
The mod~$p$ cohomology of $BPU(p)$ is detected by elementary
abelian $p$-subgroups.
\end{Thm}
\begin{proof}
According to \cite[Theorem 1.2 \& Lemma 3.1]{Moller4}, $H^*\big(PU(p)/N(PU(p))\big)$
is finite and its Euler characteristic, $\chi\big(PU(p)/N(PU(p))\big)$, equals $1$.
Therefore, the transfer argument \cite[Theorem 9.13]{DW} shows that $H^*BPU(p)\rTo
H^*BN(PU(p))$ is a monomorphism. As $H^*BN(PU(p))$ is detected by elementary abelian
subgroups by previous lemma, $H^*BPU(p)$ is so.
\end{proof}

\section{The Kono-Yagita conjecture}\label{KY-sec}

In this very short section we provide a proof of Theorem B (see Theorem \ref{demo-KY})
by means of Theorem A. Recall that for an odd prime $p$, the Milnor primitive
operators are inductively defined as $Q_0=\beta$ and $Q_{n+1}=\P{p^l} Q_n-Q_n\P{p^l}$
where $\beta$ and $\P{j}$ are the Bockstein and the $j$-th Steenrod power
respectively.

As quoted above, we use Theorem A to prove Theorem B, hence we need some information
about elementary abelian subgroups in $PU(p)$. This information is collected in the
following proposition (\cite[Corollary 3.4]{BV2} or \cite[Theorem 9.1]{AGMV})

\begin{Prop}\label{elab}
The group $PU(p)$ contains two conjugacy classes of maximal elementary abelian
subgroups corresponding to the conjugacy classes of the maximal toral elementary
abelian and a rank two nontoral.
\end{Prop}

In fact, those two subgroups of $PU(p)$ already showed up along the proof of Lemma
\ref{l2} after the identification in Lemma \ref{l1}.

The first lemma in this section shows that Conjecture \ref{KY-conj} holds for rank two
elementary abelian groups,

\begin{Lem}
Let $x$ be an odd dimensional element of $H^*B(\Z/p)^2=E(x_1,x_2)\otimes
\F_p[y_1,y_2]$, then there exists an $i$ such that $Q_m x$ is not trivial for all
$m>i$.
\end{Lem}
\begin{proof}
First notice that $Q_nx_i=y_i^{p^n}$ and $Q_ny_i=0$. Now, if $x$ is odd dimensional,
then $x=x_1f+x_2g$, where $f,g\in \F_p[y_1,y_2]$. If $Q_nx$ is nontrivial for all $n$,
lemma holds. So, let $r$ be an integer such that $Q_rx=0$. Then
$Q_rx=y_1^{p^r}f+y_2^{p^r}g=0$ and therefore there exists $h\in \F_p[y_1,y_2]$ such
that $f=y_2^{p^r}h$ and $g=-y_1^{p^r}h$. For $m>r$ we have that $$
Q_mx=y_1^{p^m}f+y_2^{p^m}g=y_1^{p^m}y_2^{p^r}h-y_2^{p^m}y_1^{p^r}h=
(y_1^{p^m-p^r}-y_2^{p^m-p^r})y_1^{p^r}y_2^{p^r}h $$ is nontrivial.
\end{proof}

Finally

\begin{Thm}\label{demo-KY}
For each odd dimensional element $x\in H^*BPU(p)$, there is $i$ such that such that
$Q_mx\ne 0$ for all $m\ge i$.
\end{Thm}
\begin{proof}
Let $x$ be in $H^*BPU(p)$ an odd dimensional element. By Theorem A, $Bj^*(x)$ is
nontrivial for some $j\colon E\rTo PU(p)$, where $E$ is an elementary abelian
$p$-group. If $E$ is toral, then $j$ factors trough maximal torus $i_T\colon T\rTo
PU(p)$. As $H^*BT$ is concentrated in even degrees, $Bj^*$ is trivial on elements of
odd degree. Therefore $Bj^*(x)$ is a non trivial odd dimensional element in $H^*BV$
for $j\colon V\rTo PU(p)$ the non toral elementary abelian subgroup which is of rank
two by Proposition \ref{elab}. By the previous lemma, there exists $i$ such that for
all $m>i$, $Q_mBj^*(x)=Bj^*(Q_mx)$ is nontrivial. Thus for all $m>i$, $Q_mx$ is
nontrivial.
\end{proof}

\section{Cohomological uniqueness}\label{demo-N}

In this section we proceed to prove Theorem C in the case $p>3$.
So in what follows $X$ is a $p$-complete space, such that there
exists an isomorphism $\phi\colon H^*BPU(p)\cong H^*X$ as an
unstable algebra over the Steenrod algebra $\mathcal{A}_p$, for
$p>3$.

The idea is to construct a homotopy equivalence $BPU(p)\pcom\rTo
X$ by means of the cohomology decomposition of $BPU(p)$ given by
$p$-stubborn subgroups \cite{JMO}.

Recall that given a compact Lie group $G$, a subgroup $P\subset G$ is called
$p$-stubborn \cite[pag.~186]{JMO} if the following conditions hold:
\begin{itemize}
\item[--] The connected component of $P$ is a torus and $\pi_0P$
is a $p$-group.
\item[--] The quotient group $N_G(P)/P$ is finite and possesses no
nontrivial normal $p$-subgroups
\end{itemize}
Then if $\mathcal{R}_p(G)$ denotes the full subcategory of the orbit category of $G$
whose objects are the homogeneous spaces $G/P$ where $P\subset G$ is $p$-stubborn, the
natural map $$ \hocolim_{G/P\in \mathcal{R}_p(G)} EG/P \rTo BG $$ induces an
isomorphism of homology with $\Z_{(p)}$-coefficients \cite[Theorem 4]{JMO}.

The $p$-stubborn subgroups of $PU(p)$ are described in the next
proposition.

\begin{Prop}\label{p-st}
The group $PU(p)$ contains exactly three $p$-stubborn subgroups up
to conjugation:
\begin{enumerate}
\item the maximal torus $T$,
\item the $p$-normalizer $N_p{\buildrel{def}\over{:=}}N_p(PU(p))$
of the maximal torus, and
\item the group $E_2=(\Z/p)^2$ generated by the diagonal matrix
$\diag(1,\zeta,\ldots, \zeta^{p-1})$, where $\zeta$ is a $p^{th}$
root of the unit, and a permutation matrix which corresponds to
the cyclic permutation $(1,2,\ldots,p)$.
\end{enumerate}
\end{Prop}
\begin{proof}
By \cite[Proposition 1.6]{JMO}, $P\subset SU(p)$ is a $p$-stubborn subgroup if and
only if $P/(P\cap Z)$ is a $p$-stubborn subgroup of $PU(p)$, where $Z\cong \Z/p$ is
the center of $SU(p)$. Finally, \cite[Theorems 6, 8 \& 10]{O} describe all the
conjugacy classes of $p$-stubborn groups in $SU(p)$, what leads to the desired result.
\end{proof}

Let $\widetilde{\mathcal{R}}_p(PU(p))$ be the full subcategory of
$\mathcal{R}_p(PU(p))$ with only the three objects: $PU(p)/T$, $PU(p)/N_p$, and
$PU(p)/E_2$. Then the strategy is to construct a homotopy commutative diagram (Lemma
\ref{HomKom}) $$\{EG/P\simeq BP\}_{PU(p)/P\in
\widetilde{\mathcal{R}}_p(PU(p))}\rTo^{f_P} X$$ such that it can be lifted to the
topological category (after Proposition \ref{obst}) so then we can recover $BPU(p)$
(up to $p$-completion) as a hocolim.

As every $p$-stubborn $P\subset PU(p)$ such that
$PU(p)/P\in\widetilde{\mathcal{R}}_p(PU(p))$ appears as a subgroup
of $N{\buildrel{def}\over{:=}}N(PU(p))$, we first construct a map
$BN\rTo X$.

\begin{Thm} \label{Mapf_N}
There exists a map $f_N\colon BN\rTo X$ such that the diagram
\begin{equation}\label{KomDgm}
\begin{diagram}
&& H^*BN&&\\ & \ruTo<{Bi_N^*}&&\luTo>{f_N^*}\\
H^*BPU(p)&&\rTo^{\phi}_\cong&&H^*X
\end{diagram}
\end{equation}
commutes.
\end{Thm}
\begin{proof}
Let $i_E\colon E=(\Z/p)^{p-1}\rTo T\rTo PU(p)$ be the maximal elementary abelian
$p$-subgroup of $PU(p)$. By Lannes' theory \cite[Th\'eor\`eme 3.1.1.]{L}, there exists
a map $f_E\colon BE\rTo X$ such that $f_E^*=Bi_E^*\phi^{-1}\colon H^*X\rTo H^*BE$. By
\cite[Proposition 3.4.6.]{L}, $$ T^E_{Bi_E^*} H^*BPU(p)\pcom\cong
H^*\map(BE,BPU(p)\pcom)_{{Bi_E}\pcom}. $$ Since $$
\map(BE,BPU(p)\pcom)_{{Bi_E}\pcom}\simeq BC_{PU(p)}(E)\pcom\simeq BT\pcom, $$ where
$C_{PU(p)}(E)$ denotes the centralizer (\cite{DZ},\cite{N3}), it follows that $$
T^E_{f_E^*} H^*X\cong T^E_{Bi^*} H^*BPU(p)\cong H^*BT\pcom. $$ Because $T^E_{f_E^*}
H^*X$ is zero in dimension $1$, we can use \cite[Th\'eor\`eme 3.2.1.]{L} and obtain $$
T^E_{f_E^*} H^*X\cong H^*\map(BE,X)_{f_E}. $$ Therefore the mapping space
$\map(BE,X)_{f_E}$ has the same cohomology ring as $BT\pcom$. The mapping space
$\map(BE,X)_{f_E}$ is $p$-complete \cite[Proposition 3.4.4]{L}, hence $BT\pcom\simeq
\map(BE,X)_{f_E}$.

Now, the standard action of $W_{PU(p)}=\Sigma_p$ on $T$ restricts to an action on $E$,
which induces an action of $\Sigma_p$ on $\map(BE,X)$. If $\sigma\in \Sigma_p$, then
$Bi_E \simeq Bi_E\sigma$, and therefore $$ f_E^*=Bi_E^*\phi^{-1}=\sigma^*
Bi_E^*\phi^{-1}=\sigma^* f_E^*, $$ and by Lannes' theory \cite[Th\'eor\`eme 3.1.1]{L},
$f_E\simeq f_E\sigma$. This means that $\Sigma_p$ acts on $\map(BE,X)_{f_E}$.

Consider now the space $Y=\map(BE,X)_{f_E}\times_{\Sigma_p} E\Sigma_p$ which fits in
the fibration $$ \map(BE,X)_{f_E}\rTo Y\rTo B\Sigma_p. $$ Fibrations with fiber
$\map(BE,X)_{f_E}$ and base $B\Sigma_p$ are classified by
$$\H^n(B\Sigma_p;\pi_n(\map(BE,X)_{f_E})=H^2(B\Sigma_p;\pi_2(\map(BE,X)_{f_E})).$$
According to \cite[Theorem 3.6]{A}, this group is trivial (recall that $p\ge 5$) which
shows that $Y\simeq BN_{p}^{\circ}$, the fiberwise $p$-completion of $BN$.

Let $f_N\colon \map(BE,X)_{f_E} \times_{\Sigma_p}E\Sigma_p\rTo X$ denote the
evaluation map. We have to prove that the diagram (\ref{KomDgm}) commutes, that is,
that $f_N^*\phi=Bi_N^*$. Let us define $a=f_N^*\phi$.

Since $H^*BN$ is detected by elementary abelian subgroups (Lemma \ref{DetElAb}), it is
enough to prove that $Bj_V^* a=Bj_V^* Bi_N^*$ for every elementary abelian $p$-group
$j_V\colon V\rTo N$. In fact, the proof of Lemma \ref{DetElAb} and Lemma \ref{l2}
shows that it is enough to consider $E$, the maximal toral elementary abelian
subgroup, and $V_n$, the nontoral elementary abelian subgroup of rank two that
coincides with $E_2$ in Proposition \ref{p-st}.

By construction of the map $f_N$, the composition $$ H^*BPU(p)\rTo^a H^*BN\rTo^{Bi^*}
H^*BT $$ is the same as $Bi_T^*$. Therefore $Bj_{E}^* a=Bj_{E}^* Bi_N^*$ for $V=E$.

Let now consider the case $V=E_2$. As $Bj_E$ factors through $BT$, $Bj_E$ can detect
only even dimensional elements in $H^*BN$. Therefore, $Bj_V$ detects $H^{odd}BN$. In
particular $Bj_V$ detects $H^3BN\supset H^3BPU(p)\subset\beta H^2BPU(p)$ which is
nontrivial (notice that $H_2(BPU(p),\Z)\cong\pi_1PU(p)=\Z/p$ so the Universal
Coefficient Theorem for cohomology \cite[Theorem 4.3 in pag.~163]{Mas} and the
description of the Bockstein morphism \cite[pag.~455]{McCleary} imply the statement).
As $f_T^*H^2X\ne 0$ by construction, then $f_N^*H^3X\subset\beta H^2X\ne 0$ and
$f_NBj_V$ detects $H^3X$ as well. Finally, by Lemma \ref{elab} the group $PU(p)$ has
only one nontoral elementary $p$-subgroup, hence by Lannes' theory there exists just
one morphism of unstable algebras $H^*BPU(p)\rTo^\psi H^*B\hat{V}$ such that $\hat{V}$
is a nontrivial elementary abelian, $H^*B\hat{V}$ is a finite module over $H^*BPU(p)$
(via $\psi$) and $\psi H^{odd}BPU(p)\ne 0$, thus $\hat{V}=V=E_2$ and
$\psi=Bj_V^*Bi_N^*$, as well as $\psi=Bj_V^* \phi f_N^*=Bj_V^* a$, thus $Bj_V^*
Bi_N^*={j'_V}^* Bi_N^*=Bj_V^* a$ also in this case.
\end{proof}

Define maps $f_P\colon EPU(p)/P\simeq BP\rTo BN\rTo^{f_N} X$ for $P=T$, $N_p$, and
$E_2$. This gives rise to a diagram
\begin{equation}\label{hom-conm}
\{EG/P\simeq BP\}_{PU(p)/P\in \widetilde{\mathcal{R}}_p(PU(p))}\rTo^{f_P} X
\end{equation}

Next lemma shows diagram \eqref{hom-conm} commutes up to homotopy.

\begin{Lem}\label{HomKom}
For every two objects $P$ and $Q$ in $\widetilde{\mathcal{R}}_p(PU(p))$ and morphism
$c_g\in \Mor(P,Q)$ the diagram
\begin{diagram}
BP & \rTo^{f_P} & X \\ \dTo<{Bc_g} && \dIgual \\ BQ & \rTo^{f_Q} & X
\end{diagram}
commutes.
\end{Lem}
\begin{proof}
Because every morphism in $\widetilde{\mathcal{R}}_p(PU(p))$ is a composition of an
automorphism and an inclusion, it is enough to prove that the diagram
\begin{diagram}
BP & \rTo^{f_P} & X \\ \dTo<{Bc_g} && \dIgual \\ BP & \rTo^{f_P} & X
\end{diagram}
commutes for every object $PU(p)/P$ in $\widetilde{\mathcal{R}}_p(PU(p))$. If $P=T$,
then the element $g$ is in the normalizer $N$, hence the diagram
\begin{diagram}
BP & \rMono& BN & \rTo^{f_N}& X \\ \dTo<{Bc_g} &&\dIgual&& \dIgual \\ BP & \rMono& BN
& \rTo^{f_N}& X
\end{diagram}
commutes.

Let $P=N_p$. Since $c_g(N_p)=N_p$, and $T$ is the connected component of $N_p$, also
$c_g(T)=T$, hence $g\in N$. Again we get a commutative diagram as in the previous
case.

Let $P=E_2$. Then $Bi_{E_2}^*=Bc_g^*Bi_{E_2}^*$, since $Bi_{E_2}\simeq Bi_{E_2}Bc_g$,
and therefore $$ f_{E_2}^*=Bi_{E_2}^*\phi^{-1}=Bc_g^* Bi_{E_2}^* \phi^{-1}=Bc_g^*
f_{E_2}^*. $$ By Lannes' theory \cite[Th\'eor\`eme 3.1.1.]{L}, $f_{E_2}\simeq f_{E_2}
Bc_g$, which finishes the proof.
\end{proof}

The diagram \eqref{hom-conm} commutes only up to homotopy, hence we do not know if the
collection of maps $\{f_P\}_{PU(p)/P\in \widetilde{\mathcal{R}}_p(PU(p))}$ induces a
map $$ \hocolim_{PU(p)/P\in \widetilde{\mathcal{R}}_p(PU(p))} EPU(p)/P\rTo X. $$ The
obstructions lie in the groups $$
\sideset{}{^i}\lim_{\overleftarrow{\widetilde{\mathcal{R}}_p(PU(p))}}
\pi_j(\map(BP,X)_{f_P}), $$ where $\lim^i$ is the $i$-th derived functor of the
inverse limit functor (\cite{BK} and \cite{WO}). Now we will prove that all
obstruction groups are trivial.

Let $$ \Pi_j^X,\Pi_j^{PU(p)}\colon \widetilde{\mathcal{R}}_p(PU(p)) \rTo\mathcal{A}b
$$ be functors defined by $$ \Pi_j^X(PU(p)/P)=\pi_j(\map(BP,X)_{f_P}), $$ $$
\Pi_j^{PU(p)}(PU(p)/P)= \pi_j(\map(BP,BPU(p)\pcom)_{(Bi_P)\pcom}), $$ where
$\mathcal{A}b$ is the category of abelian groups. Note that
$\map(BP,BPU(p)\pcom)_{(Bi_P)\pcom}\cong BZ(P)\doscom$ \cite[Theorem 3.2]{JMO}
therefore $\Pi_1^{PU(p)}(PU(p)/P)$ is well defined and, by the next lemma, also
$\Pi_1^{X}(PU(p)/P)$ is well defined.

\begin{Lem}
There exists a natural transformation $\mathcal{T}\colon \Pi_j^{PU(p)} \rTo \Pi_j^X$
which is an equivalence.
\end{Lem}
\begin{proof}
Let $P$ be the maximal torus $T$ of the $p$-normalizer $N_p$, and let
$E\cong(\Z/p)^{p-1}$ be the maximal toral elementary abelian subgroup in $N$. We apply
Lannes' $T$ functor to the diagram
\begin{equation}\label{KomDgm1}
\begin{diagram}
&& H^*BN&&\\ & \ruTo<{Bi_N^*}&&\luTo>{f_N^*}\\ H^*BPU(p)&&&&H^*X\\
\end{diagram}
\end{equation}
and get
\begin{diagram}
&& T^E_{Bi_E^*} H^*BN&&\\ & \ruTo&&\luTo&\\ T^E_{Bi_E^*} H^*BPU(p)&&&&T^E_{f_E^*}
H^*X\\
\end{diagram}
By \cite[Th\'eor\`eme 3.4.5]{L}, it follows that $$ T^E_{Bi_E^*} H^*BN\cong
H^*BC_T(E)=H^*BT, $$ $$ T^E_{Bi_E^*} H^*BPU(p)\cong H^*BC_{PU(p)}(E)=H^*BT, $$ and the
left map in the above diagram is an isomorphism. Because $T^E_{f_E^*} H^*X\cong
T^E_{Bi_E^*} H^*BPU(p)\cong H^*BT$, it is zero in the degree $1$, hence by
\cite[Th\'eor\`eme 3.2.1.]{L}, $T^E_{f_E^*} H^*X\cong H^*\map(BE,X)_{f_E}$ and the
right map in the diagram is an isomorphism. We conclude that in the diagram
\begin{diagram}
&& \map(BE,BN_p^\circ)_{(Bi_E)\pcom}&&\\ & \ldTo&&\rdTo&\\
\map(BE,BPU(p)\pcom)_{(Bi_E)\pcom}&&&&\map(BE,X)_{f_E}\\
\end{diagram}
both maps are equivariant mod~$p$ equivalences. Taking homotopy fixed points we obtain
the following diagram
\begin{diagram}
&&\map(BE,BN_p^\circ)_{(Bi_E)\pcom}^{h(P/E)}&&\\ &\ldTo&&\rdTo&\\
\map(BE,BPU(p)\pcom)_{(Bi_E)\pcom}^{h(P/E)}&&&&\map(BE,X)_{f_E}^{h(P/E)}.
\end{diagram}
where both maps are mod~$p$ equivalences, since an equivariant mod~$p$ equivalence
between 1-connected spaces induces a mod~$p$ equivalence between the homotopy
fixed-point sets. Using $\map(BP,\cdot)\simeq \map(BE,\cdot)^{h(P/E)}$, we obtain
mod~$p$ equivalences
\begin{diagram}
&&\map(BP,BN_p^\circ)_{(Bi_P)\pcom}&&\\ &\ldTo&&\rdTo&\\
\map(BP,BPU(p)\pcom)_{(Bi_P)\pcom}&&&&\map(BP,X)_{f_P}.
\end{diagram}

Let us consider the remaining case $P=E_2\cong (\Z/2)^2$. Applying Lannes' functor to
diagram (\ref{KomDgm1}) gives
\begin{diagram}
&& T^P_{Bi_P^*} H^*BN&&\\ & \ruTo&&\luTo\\ T^P_{Bi_P^*} H^*BPU(p)&&&&T^P_{f_P^*}
H^*X\\
\end{diagram}
By \cite[Th\'eor\`eme 3.4.5]{L}, we get $$ T^P_{Bi_P^*} H^*BN\cong H^*BC_N(P)=H^*BP,
$$ $$ T^P_{Bi_P^*} H^*BPU(p)\cong H^*BC_{PU(p)}(P)=H^*BP, $$ and the left map is an
isomorphism. Since $T^P_{f_P^*}H^*X$ is free in dimension $\le 2$, it follows by
\cite[Th\'eor\`eme 3.2.4]{L} that $T^P_{f_P^*} H^*X\cong H^*\map(BP,X)_{f_P}$ and also
the right map is an isomorphism. So, in the diagram
\begin{diagram}
&& \map(BP,BN_p^\circ)_{Bj_P}&&\\ & \ldTo&&\rdTo&\\ \map(BP,BPU(p)\pcom)_{Bi_P}&&
&&\map(BP,X)_{f_P}\\
\end{diagram}
both maps are mod~$p$ equivalences.

We have shown that in all cases ($P=N_p$, $T$, or $E_2$) the map $$
\map(BP,BPU(p)\pcom)_{(Bi_P)\pcom}\rTo\map(BP,X)_{f_P} $$ is a mod~$p$ equivalence.
Because $\map(BP,BPU(p)\pcom)_{(Bi_P)\pcom}$ and $\map(BP,X)_{f_P}$ are $p$-complete
spaces \cite[Proposition 3.4.4.]{L}, above map is a homotopy equivalence. To see that
this homotopy equivalence is natural, we have to prove that it commutes with maps
induced by conjugation, which means that we have to show that the diagram
\begin{diagram}
&&\map(BP,BN_p^\circ)_{(Bi_P)\pcom Bc_g}&&\\ &\ldTo&&\rdTo&\\
\map(BP,BPU(p)\pcom)_{(Bi_P)\pcom}&\lTo&\map(BP,BN_p^\circ)_{(Bi_P)\pcom}&\rTo&\map(BP,X)_{f_P}
\end{diagram}
commutes. This follows from two commutative diagrams,
\begin{diagram}
BC_{PU(p)}(P)\pcom &&&&\rIgual& BC_N(P)\pcom \\ &\rdTo^\simeq &&&&\dTo>\simeq\\
&&\map(BP,BPU(p)\pcom)_{(Bi_P)\pcom}&&\lTo &\map(BP,BN_p^\circ)_{(Bi_P)\pcom}\\
\dTo<{Bc_g}&&\dIgual&&&\\ &&\map(BP,BPU(p)\pcom)_{(Bi_P)\pcom Bc_g}&&\lTo
&\map(BP,BN_p^\circ)_{(Bi_P)\pcom Bc_g}\\ &\ruTo^\simeq &&&&\uTo>\simeq\\
BC_{PU(p)}(c_g(P))\pcom &&&&\rIgual& BC_N(c_g(P))\pcom \\
\end{diagram}
and
\begin{diagram}
BC_N(P)\pcom &\rIgual&&&& BC_{PU(p)}(P)\pcom \\ \dTo>\simeq&&&&\ldTo^\simeq &\\
\map(BP,BN_p^\circ)_{(Bi_P)\pcom} &\rTo&&\map(BP,X)_{f_P} &&\\
&&&\dIgual&&\dTo<{Bc_g}\\ \map(BP,BN_p^\circ)_{(Bi_P)\pcom Bc_g}&\rTo&&\map(BP,X)_{f_P
Bc_g} && \\ \uTo>\simeq&&&&\luTo^\simeq & \\ BC_N(c_g(P))\pcom&\rIgual&&&&
BC_{PU(p)}(c_g(P))\pcom  \\
\end{diagram}
which can be glued together.
\end{proof}

\begin{Prop}\label{obst}
For all $i,j\ge 1$,
$$\sideset{}{^i}\lim_{\overleftarrow{\widetilde{\mathcal{R}}_p(PU(p))}}
\pi_j(\map(BP,X)_{f_P})=0.$$
\end{Prop}
\begin{proof}
By the previous lemma,
$$\sideset{}{^i}\lim_{\overleftarrow{\widetilde{\mathcal{R}}_p(PU(p))}}
\pi_j(\map(BP,X)_{f_P})=
\sideset{}{^i}\lim_{\overleftarrow{\widetilde{\mathcal{R}}_p(PU(p))}}
\pi_j(\map(BP,BPU(p)\pcom)_{(Bi_P)\pcom})$$ and the right side is $0$ \cite[Theorem
4.8]{JMO}.
\end{proof}

Because all obstructions vanish, there exists a map $f\colon BPU(p)\pcom \rTo X$. By
construction of the map $f$, the diagram
\begin{diagram}
&& (BN_p)\pcom&&\\ & \ldTo<{Bi_N}&&\rdTo>{f_N}\\ BPU(p)\pcom&&\rTo^{f}&&X\\
\end{diagram}
commutes. The Euler characteristic $\chi(PU(p)/N_p)\ne 0 \mod p$, hence a transfer
argument shows that $Bi_N^*$ is a monomorphism. By Theorem \ref{Mapf_N}, also $f_N^*$
is a monomorphism. Therefore, $f^*$ is a monomorphism and, because $H^*BPU(p)\cong
H^*X$ is a finite dimensional in each degree, $f^*$ is an isomorphism. This shows that
$f$ is a homotopy equivalence.


\end{document}